\documentclass[11pt,a4]{amsart}
\usepackage{amsmath,amssymb}
\makeatletter
\renewcommand{\section}{\@startsection{section}{1}{0pt}{20pt}{6pt}{\large\bfseries}}
\makeatother

\usepackage{enumerate}

\numberwithin{equation}{section}

\theoremstyle{plain}
  \newtheorem{thm}{Theorem}[section]
  \newtheorem{prop}[thm]{Proposition}
  \newtheorem{lemma}[thm]{Lemma}
   \newtheorem{cor}[thm]{Corollary}
\theoremstyle{definition}
  \newtheorem{remark}[thm]{Remark}

\newcommand{\R}{\mathfrak{R}}
\renewcommand{\Re}{{\mathfrak{Re}}}
\newcommand{\C}{\mathfrak{C}}
\newcommand{\Q}{\mathbb{Q}}
\newcommand{\E}{\mathbb{E}}

\renewcommand{\L}{\mathbf{L}}

\renewcommand{\P}{{\mathbb{P}}}

\renewcommand{\d}{{\rm{d}}}
\newcommand{\Ex}{\Sigma}
\newcommand{\iex}{A}

\newcommand{\M}{\mathcal{E}}
\newcommand{\K}{\mathcal{LK}}

\newcommand{\I}{\mathcal{I}}
\newcommand{\Ip}{\mathcal{I}_{\alpha,\psi}}
\newcommand{\Ia}{\I_{\alpha,\psi_{\rho}}}
\newcommand{\Ig}{\I_{\alpha,\psi_{\gamma}}}

\newcommand{\It}{\I_{\alpha,\psi_{\theta}}}

\newcommand{\an}{a_n(\psi_{\gamma};\alpha)}
\newcommand{\ann}{a_n(\psi;\alpha)}

\newcommand{\N}{\mathcal{N}}

\newcommand{\Id}[1]{{{\mathbb{I}}}_{\{#1\}}}

\begin{document}
\title{Infinite divisibility of solutions to some self-similar integro-differential equations and
 exponential functionals of L\'evy processes}
\author{P. Patie}
\address{Institute of Mathematical Statistics and Actuarial Science, University of Bern,
Alpeneggstrasse, 22, CH-3012 Bern, Switzerland}
\email{patie@stat.unibe.ch}
\thanks{I wish to thank an anonymous referee for very valuable and helpful comments  that led to improving the presentation of the paper. I am also grateful to V. Rivero for many fruitful discussions on the topic during my visit of CIMAT in Guanajuato.}
\begin{abstract}
We first characterize the increasing eigenfunctions associated to
the following family of integro-differential operators, for any
$\alpha,x >0, \gamma \geq 0$ and $f$ a smooth function on $\R^+$,
\begin{eqnarray} \label{eq:inf_aba}
\mathbf{L}^{(\gamma)} f(x)&=& x^{-\alpha}
\left(\frac{\sigma}{2} x^{2}f''(x) + (\sigma\gamma +b) x
f'(x)\right.
\\&+&
\left.\int^{\infty}_{0}\left((f(e^{-r}x)-f(x))e^{-r\gamma}+xf'(x)r\Id{r\leq1}\right)\nu(dr)\right)
\nonumber
\end{eqnarray}
where the coefficients $ b \in \R, \: \sigma \geq 0$ and the measure
$\nu$, which satisfies the integrability condition $\int^{\infty}_0
(1\wedge r^2 )\: \nu(dr) <+ \infty$, are uniquely determined by the
distribution of a spectrally negative infinitely divisible random
variable, with characteristic exponent $\psi$.
$\mathbf{L}^{(\gamma)}$ is known to be the infinitesimal generator
of a positive $\alpha$-self-similar Feller process,
which has been introduced by Lamperti \cite{Lamperti-72}. The eigenfunctions are
expressed in terms of a new family of power series which includes,
for instance, the modified Bessel functions of the first kind and
some generalizations of the Mittag-Leffler function. Then, we show
that some specific combinations of these functions are Laplace
transforms of self-decomposable or infinitely divisible
distributions concentrated on the positive line with respect to the main argument, and, more surprisingly, with respect to the parameter
$\psi(\gamma)$. In particular, this generalizes a result of
Hartman \cite{Hartman-76} for the increasing solution of the Bessel differential
equation. Finally,  we compute, for some cases, the
associated decreasing eigenfunctions and derive the Laplace
transform of the exponential functionals of some spectrally negative
L\'evy processes with a negative first moment.
\end{abstract}


\maketitle

\bibliographystyle{amsplain}

\section{Introduction}
During the last decade, there has been a  renewed interest for
self-similar semigroups, something which seems to be attributed to
their
 connections to several fields
of mathematics and more generally to many area of sciences. For
instance, in probability theory, these semigroups arise in the study
of important processes such as self-similar processes, branching
processes and also in the investigation of self-similar
fragmentation. Moreover, the Feller processes associated to
self-similar semigroups are closely related, via the Lamperti's
mapping, to the exponential functionals of L\'evy processes which
appear to be  key objects in a variety of settings (random processes
in random environment, mathematical finance, astrophysics $\ldots$).
We refer to Bertoin and Yor \cite{Bertoin-Yor-05} for an interesting recent
 survey on this topic. Finally, we emphasize that they are
also related with the theory of fractional operator which is used
intensively in many applied fields, see e.g.~the survey paper of
Kilbas and Trujillo \cite{Kilbas-Trujilo-01}. \\In this paper, we
provide, in terms of power series, the increasing eigenfunction, associated to the linear operator
$\mathbf{L}^{(\gamma)}$, given by \eqref{eq:inf_aba},
that is the increasing solution to the integro-differential
equation, for $x,q\geq0$,
\begin{equation*}
\mathbf{L}^{(\gamma)} f_q(x) = q f_q(x).
\end{equation*}
As a byproduct, we compute the Laplace transform of the first
passage times above for spectrally negative self-similar processes
and some related quantities. \\
Moreover, when the spectrally negative random variable has a negative first moment, we
provide, under an additional technical condition, the decreasing eigenfunctions associated to
$\mathbf{L}=\mathbf{L}^{(0)}$. We deduce an explicit
expression of the Laplace transform of the exponential
functional of some spectrally negative L\'evy processes with negative mean. This is a companion result of Bertoin and Yor \cite{Bertoin-Yor-02} who characterized, in terms of its negative entire moments, the law of the exponential
functional of spectrally positive L\'evy processes which drift to $-\infty$.

Furthermore, it is plain that $\mathbf{L}^{(\gamma)}$ is a generalization of the
infinitesimal generator of  the (re-scaled) Bessel process
($\nu\equiv 0$ and $\alpha=2$). In this specific case, Hartman
\cite{Hartman-76}, relying on purely analytical arguments, showed, that
the function
\begin{eqnarray} \label{eq:hr}
  \gamma \mapsto
  \frac{{\rm{I}_{\sqrt{2\gamma}}}(a){\rm{I}_0}(A)}{{\rm{I}_{\sqrt{2\gamma}}}(A){\rm{I}_0}(a)},
  \quad 0<a<A<\infty,
\end{eqnarray}
is the Laplace transform of an infinitely divisible distribution concentrated on
the positive line, where $\rm{I}_{\nu}$ stands for the modified Bessel function
of the first kind.  We mention that in the limit case $A \rightarrow
\infty$, the result above has been reproved, in an elegant fashion,
by Pitman and Yor \cite{Pitman-Yor-80}. We shall provide a simple probabilistic explanation of Hartman's result (for any $0<A<\infty$) and show that this property still hold for similar ratios of the increasing eigenfunctions associated to the non-local operator $\mathbf{L}^{(\gamma)}$.

The outline of the remaining of the paper is as follows. In the sequel, we set up the notation and provide some basic results.
Section 2 is devoted to the statement of the main results. The
proofs are given in Section 3. Finally, in the last section, we
illustrate our approach by investigating some known and new
examples.

\subsection{Notation and preliminaries}
\subsubsection{Some important set of probability measures}  Let $\xi_1$ be a spectrally negative infinitely divisible
random variable. It is well known that
its characteristic exponent, $\psi$, admits the following
L\'evy-Khintchine representation
\begin{eqnarray}\label{eq:lap-levy}
\psi(u) = b u + \frac{\sigma}{2} u^2 + \int^{\infty}_0 (e^{-u r} -1
+u r\Id{r\leq1})\nu(dr), u\geq0,
\end{eqnarray}
where the coefficients $ b \in \R, \: \sigma \geq 0$ and the measure
$\tilde{\nu}$, image of $\nu$ by the mapping $x \rightarrow -x$,
which satisfies the integrability condition $\int_{-\infty}^0
(1\wedge r^2 )\: \tilde{\nu}(dr) <+ \infty$, are uniquely
determined by the distribution of $\xi_1$. We exclude the case when $b\leq 0$ and $\int_{-\infty}^0
(1\wedge r )\: \tilde{\nu}(dr) <+ \infty$, i.e. when $\psi$ is the Laplace exponent of the negative of a subordinator. It is plain that $\lim_{u \rightarrow
\infty}\psi(u)=+\infty$ and by monotone convergence, one gets
$\E[\xi_1]= b - \int_1^{\infty}r\nu(dr) \in [-\infty,\infty).$
Differentiating again, one observes that $\psi$ is strictly convex
unless $\xi$ is degenerate, which we exclude. Note that $0$ is
always a root of the equation $\psi(u)=0$. However, in the case
$\E[\xi_1] <0$, this equation admits another positive root, which we denote
by $\theta$. This yields the so-called Cram\'er condition
\begin{eqnarray*}
\E[e^{\theta \xi_1}] = 1.
\end{eqnarray*}
Then, for any $\E[\xi_1]  \in [-\infty,\infty)$, the function $u\mapsto \psi(u)$ is continuous and increasing on $[\max(\theta,0),\infty)$ and thus it has a
well-defined inverse function $\phi:[0,\infty)\rightarrow
[\max(\theta,0),\infty)$ which is also continuous and  increasing. We
denote the totality of all functions $\psi$ of the form
\eqref{eq:lap-levy} by $\K$. Note, from the stability of infinitely
divisible distributions under convolution and convolution powers to
positive real numbers, that $\K$ forms a convex cone in the space of
real valued functions defined on $[0,\infty)$.

We also mention that when a probability measure $dm$ is supported on a subset of $[0, \infty)$, then $dm$ is infinitely divisible if and only if
its Laplace transform satisfies the conditions, for any $u\geq0$,
\begin{equation*}
\int_0^{\infty}
e^{-u x} dm(x)= e^{- \phi(u)}
\end{equation*}
with $\phi(0)= 0$ and $\phi'(u)$ is completely monotonic, i.e.~$\phi'$  is infinitely differentiable  on $(0,\infty)$ and for all $n=1,2\ldots$, $(-1)^{n-1}\phi^{(n)}(u)> 0,\: u > 0$.
Moreover, the so-called Laplace exponent, $\phi$ admits the following L\'evy-Khintchine representation
\begin{eqnarray} \label{eq:sub}
\phi(u) = au +  \int_0^{\infty} (1-e^{-u r})\mu(dr), \: u\geq0,
\end{eqnarray}
for some $a\geq0$ and some positive measure $\mu$ on $(0,\infty)$ satisfying $\int_0^{\infty}
(1\wedge r)\mu(dr) <\infty$, see e.g.~Meyer \cite{Meyer-69}.

Finally, we recall that a random variable $H$ is
self-decomposable (or of class $L$) if it is solution to the random
affine equation
\begin{equation*}
H \stackrel{(d)}{=} cH + H_c
\end{equation*}
where  $\stackrel{(d)}{=}$ stands for the equality in distribution,
$0<c<1$  and $H_c$ is a random variable
independent of $H$.  It is well-known
that the law of these random variables is absolutely continuous, see e.g.~\cite[Example 27.8]{Sato-99}.
Moreover, Wolfe \cite{Wolfe-71} showed that the density $h$ of a
positive self-decomposable random variable is unimodal, i.e.~there
exists $a \in \R^+$ (the mode) such that  $h$ is increasing on $]0,
a[$ and decreasing on $]a,\infty[$. The Laplace exponent,
$\phi_s$, of a self-decomposable distribution concentrated on $\R^+$
is given by
\begin{eqnarray*}
\phi_s(u) = au +  \int_0^{\infty} (1-e^{-u r})\frac{k(r)}{r}dr
\end{eqnarray*}
where $k$ is a positive decreasing function.
\noindent We refer to the monographs of Sato \cite{Sato-99} and Steutel and van Harn \cite{Steutel-vanHarn-04} for an
excellent account on these sets of probability measures. \\

\subsubsection{L\'evy and Lamperti}
Let $\P_x$ (we write simply $\P$ for $\P_0$) be the law  of a spectrally negative L\'evy process
$\xi:=(\xi_t)_{t\geq 0}$, starting at $x \in \R$, with
$(F_t)_{t\geq0}$ its natural filtration. This law is characterized
by the characteristic exponent of $\xi_1$, which we assumed to
belong to $\K$, i.e.~being of the form \eqref{eq:lap-levy}. We
deduce, from the above discussion and the strong law of large
numbers, that $\lim_{t \rightarrow +\infty} \xi_t =
\textrm{sgn}(\E[\xi_1])\infty$ a.s.~and the process oscillates if $\E[\xi_1]=0$.\\
\noindent For any $\gamma \geq 0$, we write $\P^{(\gamma)}$ for the law
of the L\'evy process with characteristic  exponent
\begin{equation*}
\psi_{\gamma}(u) = \psi(u+\gamma)-\psi(\gamma), \: u\geq0.
\end{equation*} The
laws $\P^{(\gamma)}$ and $\P$ are connected via the following
absolute continuity relationship, also known as  Esscher transform,
\begin{equation} \label{eq:esscher}
d\P_x^{(\gamma)}{}_{\mid F_t}=e^{\gamma
(\xi_t-x)-\psi(\gamma)t}d\P_x{}_{\mid F_t}, \quad t>0,\:x\in \R.
\end{equation}
\noindent Lamperti \cite{Lamperti-72} showed that there exists a one to one
mapping between $\P_x$ and the law
$\Q_{e^{x}}$ of a $\frac{1}{\alpha}$-self-similar Markov process
$X$ on $(0,\infty)$, i.e.~a Feller process which enjoys the
following $\alpha$-self-similarity property, for any $c>0$,
\begin{equation} \label{eq:self}
\left((X_{c^{\alpha}t})_{t\geq0}, \Q^{(\gamma)}_{ce^{x}}\right)
\stackrel{(d)}{=}\left((cX_{t})_{t\geq0},
\Q^{(\gamma)}_{e^x}\right).
\end{equation}
More precisely, Lamperti showed that $X$ can be constructed from
$\xi$ as follows
\begin{equation} \label{eq:lamp_transf}
\log\left(X_t\right) =  \xi_{\iex_t}, \: t\geq 0,
\end{equation}
 where
\begin{equation*}
 \iex_t = \inf
\{ s \geq 0; \: \Ex_s := \int_0^s e^{\alpha \xi_u} \: du
> t \}.
\end{equation*}  We write ${\rm{E}}^{(\gamma)}_x$ (resp.~${\rm E}_x$) the
expectation operator associated with $\Q^{(\gamma)}_x$
(resp.~$\Q_x=\Q_x^{(0)}$).  Moreover, for $\E[\xi_1] <0$,  it is
plain that $X$ has a a.s.~finite lifetime which is
$\kappa_0=\inf\{s\geq0;\: X_{s^-}=0, \rm{ } X_s=0\}$. However, under the additional
condition $0<\theta<\alpha$, where we recall that $\psi(\theta)=0$,
Rivero \cite{Rivero-05}, showed that the minimal process
$(X,\kappa_0)$ admits a \emph {unique recurrent extension that hits and
leaves $0$ continuously a.s.~}and which is a
$\alpha$-self-similar process on $[0,\infty)$. We simply
write $(X,\Q_x)$ for the law of such a recurrent extension starting from
$x\geq0$. Furthermore, for $\E[\xi_1]\geq 0$, Bertoin and Yor
\cite[Proposition 1]{Bertoin-Yor-02}, showed that the family of
probability measures $(\Q_x)_{x>0}$ converges in the sense of finite
dimensional distribution to a probability measure $\Q_{0^+}$ as $x
\rightarrow 0+$, see also Caballero and Chaumont \cite{Caballero-Chaumont-06} for conditions for  the weak convergence. Thus, for any $x\geq0$, $(X,\Q_x)$ is also
spectrally negative, in the sense that it has no positive jumps.
Moreover,  for any $x\geq0$, $(X,\Q_x)$ is a Feller process on
$[0,\infty)$ and we denote its semigroup (resp.~its resolvent) by
$(Q_t)_{t\geq0}$ (resp.~by $U^q, q>0$), i.e.~for any $x,t \geq0$ and $v \in
B([0,\infty))$, the space of bounded Borelian functions on
$[0,\infty)$,
\begin{eqnarray*}
Q_t v(x) &=& {\rm E}_x[v(X_t)] \\
U^q v(x) &=& \int_0^{\infty}e^{-qt}Q_t v(x)dt.
\end{eqnarray*}
We also introduce the semigroup and the resolvent of the minimal process $(X,\kappa_0)$, for $x>0$,
\begin{eqnarray*}
Q^0_t v(x) &=& {\rm E}_x[v(X_t), t<\kappa_0] \\
U_0^q v(x) &=& \int_0^{\infty}e^{-qt}Q^0_t v(x)dt.
\end{eqnarray*}
The strong Markov property yields the following expression for the resolvent of the recurrent extension
\begin{equation} \label{eq:res-e}
U^q v(x) = U_0^qv(x) +\E_x[e^{-q\kappa_0}] U^q v(0), \: x\geq0.
\end{equation}

Moreover, we recall that the Lamperti mapping reads in terms of
the characteristic operator  as follows.
\begin{prop}
Let $f : \R^+ \rightarrow \R$ be such that $f(x),xf'(x)$ and
$x^2f''(x)$ are continuous functions on $\R^+$, then $f$ belongs to
the domain, $\mathbb{D}(\L)$, of the characteristic operator $\L$ of $(X,\Q)$ which is given, for $x>0$, by
\begin{eqnarray}  \label{eq:inf_ab}
\mathbf{L} f(x) &=& x^{-\alpha}
\left(\frac{\sigma}{2} x^{2}f''(x) + b x
f'(x)\right.
\\&+&
\left.\int^{\infty}_{0}\left((f(e^{-r}x)-f(x))+xf'(x)r\right)\nu(dr)\right).
\nonumber
\end{eqnarray}
In the case $\E[\xi_1] <0$,  $f \in \mathbb{D}(\mathbf{L})$ if and only if it satisfies the boundary condition
\begin{equation} \label{eq:cond-b}
\lim_{x \rightarrow 0} \frac{f'(x)}{x^{\theta-1}} =0.
\end{equation}
\end{prop}
\begin{proof}
The first part of the proposition follows from Lamperti
\cite[Theorem 6.1]{Lamperti-72}. We point out that, in the case
$\E[\xi_1] <0$, the characteristic operator of the minimal process
and its recurrent extension coincide for $x>0$. Indeed, from its
Feller property, the semigroup $(Q_t, t>0)$ corresponding to the
recurrent extension leaves invariant $C_0([0,\infty))$, the space of
continuous functions vanishing at infinity. Therefore, if we let
\eqref{eq:inf_ab} stand for the characteristic operator of
the process, then we can conclude, see Gikhman and Skorohod \cite[p.~130, Theorem
1]{Gikhman-Skorohod-75}, that the domain $D(L)$ of the strong
infinitesimal generator, $L$, of the process $(X,\Q)$ consists of all the
functions $f \in C_0([0,\infty))\cap
\mathbb{D}(\mathbf{L})$ such that $\mathbf{L}f \in
C_0([0,\infty))$.  However, it is plain that the
analytic form of the function $\mathbf{L}f(x)$ for $x > 0$ and for
any function $f$, such that $f(x),xf'(x),x^2f''(x)$ are continuous, does not
depend on the method of extension and is given by the expression
above.
Let us now turn to the boundary condition in the case $\E[\xi_1] <0$. We first deal with the necessary condition. From the discussion above, it is clear that it is enough to characterize the boundary condition for the strong infinitesimal generator. To this end, we recall that $U^q$ is a Fellerian resolvent, see Rivero \cite[Theorem
2]{Rivero-05}. Thus, let $f \in C_{0}([0,\infty)) \cap
\mathbb{D}(\mathbf{L})$ then there exists $v\in C_{0}([0,\infty)) $
such that $\mathbf{L}f-qf=v$, i.e.~for $x\geq0$
\begin{equation*}
U^q v(x) = f(x).
\end{equation*}
Then, using the expression of the resolvent of the recurrent extension \eqref{eq:res-e} and the continuity of $f$, we get
\begin{equation} \label{eq:der}
f(x) - f(0) = U_0^qv(x) -\E_x[1-e^{-q\kappa_0}] U^q v(0).
\end{equation}
Let us now consider $v_0 \in C_{0}((0,\infty))$, the space of  continuous function vanishing at $0$ and $\infty$. Then, from e.g.~\cite[Lemma 1]{Rivero-05}, we have $U^qv_0(x) \in C_{0}((0,\infty)) \cap \mathbb{D}(\mathbf{L})$.
Moreover, from  \cite[Theorem 2]{Rivero-05}, we know that a necessary condition  for the existence of a unique recurrent extension which hits and leaves $0$ continuously a.s.~is that both limits
\begin{equation*}
\lim_{x \rightarrow 0} \frac{U_0^qv_0(x)}{x^{\theta}} \quad {\rm{ and }} \quad \lim_{x \rightarrow 0} \frac{\E_x[1-e^{-q\kappa_0}]}{x^{\theta}}
\end{equation*}
exist for any $v_0 \in C_0((0,\infty))$. Moreover, in this case, the identity
\begin{equation} \label{eq:el}
\lim_{x \rightarrow 0} \frac{U_0^qv_0(x)}{\E_x[1-e^{-q\kappa_0}]}= U^q v_0(0)
\end{equation}
holds for any $v_0 \in C_0((0,\infty))$.
Hence, for $v_0 \in C_0((0,\infty))$,  the necessary condition  \eqref{eq:cond-b} follows  by dividing both sides of \eqref{eq:der} by $x^{\theta}$, by taking the limit $x \rightarrow 0$ and by invoking the uniqueness of the limit. The general case, i.e.~for any $v\in C_{0}([0,\infty)) $, follows from the fact that the Lebesgue measure of the set $\{t\geq0;X_t=0\}$ is, by construction,  $0$ with probability $1$, and from Blumenthal \cite[Section 4]{Blumenthal-83}. The sufficient part is readily obtained from the uniqueness of the recurrent extension that hits and
leaves $0$ continuously a.s.
\end{proof}

\subsubsection{The family of power series} Let $\psi \in \K$ and for $\gamma \geq 0$ and $\alpha>0$,
set
\begin{equation*}
\an=\left(\prod_{k=1}^n \psi_{\gamma}(\alpha k)\right)^{-1}, \quad a_0=1,
\end{equation*}
where we recall that $\psi_{\gamma}(u) = \psi(u+\gamma)-\psi(\gamma), \: u\geq0$.
Then, we introduce the function $\Ig$ which admits the series
representation
\begin{equation*}
\Ig(z)=\sum_{n=0}^{\infty} \an z^{n}, \quad
 z
\in  \C.
\end{equation*}
We simply write $\Ip$ when $\gamma=0$. We gather some basic properties of this family  of power series
which will be useful for the sequel.
\begin{prop} \label{prop:pow}
For
any $\psi \in \K$, $\gamma \in \C$, $\Re(\gamma) \geq 0$ and
$\alpha >0$, $\Ig$ is an entire function. Moreover, if $\E[\xi_1]\geq 0$ or if $\E[\xi_1]<0$ and $\theta
<\alpha$, $\Ip$ is positive and increasing on
$[0,\infty)$.
\end{prop}
\begin{proof}
Observe that
\begin{eqnarray*}
\frac{\mid a_{n+1} \mid}{\mid a_{n} \mid}&=& \frac{1}{\mid
\psi(\alpha n+\gamma)-\psi(\gamma)\mid}.
\end{eqnarray*}
The analyticity of $\Ig$ follows from the fact that $\lim_{u
\rightarrow \infty}\psi_{\gamma}(u) = +\infty$. The positivity and the
monotonicity is secured by observing that under the condition
$\theta <\alpha$ in the case $\E[\xi_1] <0$, we have $\psi(\alpha)>0$ and  $\psi$ is increasing on $[\max(\theta,0),\infty)$.
\end{proof}
\begin{remark}
Note that Rivero's condition, $0<\theta<\alpha$ for $\E[\xi_1] <0$, arises
naturally in the previous proposition to ensure that the associated
functions are positive and increasing.
\end{remark}

\section{Main results}
Let $\psi \in \K$. Moreover, if  $\E[\xi_1] <0$,  we assume that
$\theta<\alpha$,  recalling that $\psi(\theta)=0.$ Next, for
 $a \in \R$, we introduce the stopping times
\begin{eqnarray*}
\tau_a =\inf\{s \geq 0; \: \xi_s = a\} \textrm{ and } \kappa_{e^a}
=\inf\{s \geq 0; \: X_s = e^a\}
\end{eqnarray*}
with the convention that $\inf\{\O\}= \infty$. For any $\lambda \geq
0$, we denote $\rho=\phi(\lambda)$ where $\phi:[0,\infty)\rightarrow
[\max(\theta,0),\infty)$  is the increasing and continuous inverse
function of $\psi$.

\begin{thm} \label{thm:1} Let $q \geq 0$ and $0\leq x\leq a$. Then, we have
\begin{eqnarray} \label{eq:lap_ssm}
 {\rm{E}}_x \left[e^{-q \kappa_a } \right]
 &=&\frac{\Ip(qx^{\alpha})}{\Ip(qa^{\alpha})}.
\end{eqnarray}
Moreover, for $\lambda\geq 0$,
\begin{equation} \label{eq:jla}
{\rm{E}}_x\left[e^{-q\kappa_a-\lambda\iex_{\kappa_a}}\Id{\kappa_a<\kappa_0}\right]=
\left(\frac{x}{a}\right)^{\rho}\frac{\Ia\left(qx^{\alpha}\right)}{\Ia\left(qa^{\alpha}\right)}
\end{equation}
and
\begin{equation} \label{eq:jlp}
\E_x \left[e^{-\lambda
\tau_a-q\Sigma_{\tau_a}}\Id{\tau_a<+\infty}\right] =
e^{\rho(x-a)}\frac{\Ia\left(qe^{\alpha
x}\right)}{\Ia\left(qe^{\alpha a}\right)}.
\end{equation}
\end{thm}
\begin{remark}
\begin{enumerate}
\item
Note, by letting $q \rightarrow 0$ in \eqref{eq:lap_ssm}, that
$\Q_x\left[\kappa_a<+\infty\right]=1$ for any $0\leq x \leq a$.
Thus, the points above the starting point are recurrent states for
$(X,\Q)$. Moreover, for $\E[\xi_1]\geq 0$ and $x<a$, we recall that
$\P_x\left[\tau_a<+\infty\right]=1$ and
$\Q_{e^x}\left[\kappa_{e^a}<\kappa_0\right]=1$. Thus, under such a
condition, the indicator  functions in \eqref{eq:jla} and in
\eqref{eq:jlp} can be omitted.
\item Consider the case $\psi(u)=\frac{1}{2}u^2+ \gamma u$, $\gamma>-1$, i.e.~$(X,\Q)$ is a
Bessel process of index $\gamma$. Then, the expression of the Laplace transform of $\kappa_a$, which is well-known to be expressed in terms of the modified Bessel functions of the first kind, dates back to Ciesielski and Taylor \cite{Ciesielski-Taylor-62} and Kent \cite{Kent-78}. We refer to  Section \ref{sec:b} for more
detailed computations related to the modified Bessel functions.

\end{enumerate}
\end{remark}
We proceed by characterizing, through its Laplace transform,  the law of the exponential functional, $\Sigma_{\infty}$, of spectrally negative L\'evy processes satisfying Rivero's condition. We emphasize that it is a companion result of Bertoin and Yor \cite[Proposition 2]{Bertoin-Yor-02} who computed the negative entire moments of the exponential functional of spectrally positive L\'evy  processes when they drift to $-\infty$.

\begin{thm} \label{thm:2}
Assume $\E[\xi_1] <0$ and $0<\theta<\alpha$. Then, there exists a
positive constant $C_{\theta}$ such that
\begin{equation*}
\Ip\left(x^{\alpha}\right)\sim C_{\theta}
 x^{\theta}\It\left(x^{\alpha}\right) \quad {\textrm{ as  }} \: x\rightarrow
 \infty.
 \end{equation*}
Moreover, if we assume that
there exists $\beta \in [0,1]$ such that $\lim_{u \rightarrow
\infty} \psi(u)/u^{1+\beta}= l_{\beta}$ then we have
\begin{eqnarray*}
 C_{\theta} =  \frac{\Gamma(1-\frac{\theta}{\alpha})}{\alpha
}l^{-\theta_{\alpha}}_{\beta}e^{M_{\gamma} \beta
\theta_{\alpha}}\prod_{k=1}^{\infty}e^{-\frac{\beta\theta_{\alpha}}{k}}\frac{(k+\theta_{\alpha})\psi(\alpha
k)}{k\psi(\alpha k+\theta_{\alpha})}
\end{eqnarray*}
where $M_{\gamma}=0.577\ldots$ stands for Euler-Mascheroni constant and
$\theta_{\alpha}=\frac{\theta}{\alpha} <1$. \\Next, introduce the
function
\begin{equation*}
\N_{\alpha,\psi,\theta}(x) =\Ip\left(x\right)-
C_{\theta}x^{\frac{\theta}{\alpha}}\It\left(x\right),\quad x\geq0.
\end{equation*}
$\N_{\alpha,\psi,\theta}$ is analytical on the right half plane and
 decreasing on $\R^+$. Finally, the positive random variable $\Ex_{\infty}$ has
the following Laplace transform
\begin{eqnarray} \label{eq:lt_exp}
\E \left[e^{-q \Ex_{\infty}} \right] &=&\N_{\alpha,\psi,\theta}(q).
\end{eqnarray}
\end{thm}
\begin{remark}
Note that the random variable  $\Sigma_{\infty}$ is solution to the random equation, for any $a>0$,
\begin{eqnarray*}
\Ex_{\infty} &\stackrel{(d)}{=}&  \Ex_{\tau_{-a}} + e^{-\xi_{\tau_{-a}}} \Ex'_{\infty}
\end{eqnarray*}
where $\Ex'_{\infty}$ is an independent copy of $\Ex_{\infty}$.
Indeed, together, the strong Markov property, the stationarity and
independency of the increments of the L\'evy process $\xi$ entails
that, for any $a \in \R$, the shifted process $(\xi_{t+\tau_a}- \xi_{\tau_a})_{t\geq0}$
is distributed as $(\xi_t)_{t\geq0}$ and is independent of $(\xi_t,
t\leq \tau_a)$. Finally we get the equation by noting that, since $\E[\xi_1] <0$, we have for any $a>0$, $\P(\tau_{-a} <\infty)=1$.
\end{remark}

Bertoin and Yor  \cite{Bertoin-Yor-02}  determined, in terms of their positive entire moments, the entrance law of spectrally negative self-similar positive Markov processes when $\E[\xi_1]\geq0$. In the sequel, we  characterize the entrance law of the dual
process of $(X,\Q)$ when $-\infty<\E[\xi_1]<0$, i.e.~of spectrally positive self-similar positive Markov processes when the underlying L\'evy process, in the Lamperti mapping, has a finite positive mean.  To this end, let $(\widehat{X},\Q)$ be the
self-similar process associated to the L\'evy process
$(\widehat{\xi},\P)$, the dual of $(\xi,\P)$ with respect to the
Lebesgue measure. We recall that for $-\infty<\E[\xi_1] <0$, Bertoin and Yor
\cite[Lemma 2]{Bertoin-Yor-02-b} showed that, for $x>0$,
$(\widehat{X},\Q_x)$ is weak duality, with respect to the reference
measure $m(dy)=\alpha y^{\alpha-1}dy$, with the minimal process
$(X,\kappa_0)$. In the same vein, Rivero \cite[Lemma 7]{Rivero-05}
proved that in the case $\E[\xi_1]<0$ and $\theta< \alpha$, $(X,\Q)$
is in weak duality, with respect to the measure
$m^{\theta}(dy)=y^{\alpha-\theta-1}dy$, with
$(\widehat{X},\Q^{(\theta)})$, the unique recurrent extension which
hits and leaves $0$ continuously a.s., of the self-similar process
associated, via the Lamperti's mapping, to
$(\widehat{\xi^{\theta}},\P)$, the dual of the $\theta$-Esscher
transform of $(\xi,\P)$. Before stating the next result, we recall
that an entrance law $\{\eta_s; s>0\}$ for the semi-group $Q_t$ is a
family of finite measures on the Borel sets of $(0,\infty)$ such
that $\eta_sQ_t=\eta_{s+t}$ for all strictly positive $s$ and $t$
and such that ${\rm E}^{\eta_s}[1-e^{-\kappa_0}]$ remains bounded as
$s$ approaches $0$.

\begin{cor} \label{cor:1}
If $-\infty<\E[\xi_1] <0$, then  $(\widehat{X},\Q)$ admits an
entrance law which is absolutely continuous with respect to the
reference measure $m(dy)$. Its Laplace transform with respect to the time variable is given, for
$y,q>0$, by
\begin{equation*}
\widehat{n}^q(y) = \frac{1}{|\E[\xi_1]|}
\N_{\alpha,\psi,\theta}(qy^{\alpha}).
\end{equation*}
Moreover, assume $\E[\xi_1] <0$ and $0<\theta<\alpha$. Then, $(\widehat{X},\Q^{(\theta)})$
admits an entrance law which is absolutely continuous with respect
to the reference measure $m^{\theta}(dy)$. Its Laplace transform with respect to the time variable  is
given, for $y,q>0$, by
\begin{eqnarray*}
\widehat{n}^{q}_{\theta}(y)&=&\frac{1}{\psi'(\theta) C_{\theta}}\N_{\alpha,\psi,\theta}(qy^{\alpha}).
\end{eqnarray*}
\end{cor}

Finally, we show that some specific combinations of the
functions $\Ip$ define some mappings from the convex cone $\K$ into
the convex cone of positive self-decomposable distributions or into
the convex cone of positive infinitely divisible distributions.
\begin{thm} \label{thm:3}
Let $q \geq0$. Then, the mapping
\begin{eqnarray} \label{eq:selfd}
 q \mapsto \frac{1}{\Ip(q)}
\end{eqnarray}
is the Laplace transform of a positive self-decomposable
distribution. \\
The mappings
\begin{eqnarray} \label{eq:id}
 q \mapsto  \exp\left(-q\frac{\frac{d}{d q}\Ip(q)}{\Ip(q)}\right)  \: \textrm{ and }  q \mapsto
\frac{1}{\Ip(q)} \exp\left(-q\frac{\frac{d}{dq}\Ip(q)}{\Ip(q)}\right)
\end{eqnarray}
are the Laplace transforms of positive infinitely divisible
distributions.
 \\Finally, for $\E[\xi_1]\geq 0$, $ \lambda
\geq0$, $0<a<A<\infty$ and recall that $\rho=\phi(\lambda)$, the
mapping
\begin{eqnarray} \label{eq:hw}
  \lambda \mapsto \left(\frac{a}{A}\right)^{\rho}  \frac{\Ia(a)\Ip{(A)}}{\Ia(A)\Ip(a)}
\end{eqnarray}
is the Laplace transform of an infinitely divisible distribution on
the positive line.
\end{thm}
\begin{remark} \label{rem:hw}
The random variable which Laplace transform is given in \eqref{eq:hw} is
characterized below in \eqref{eq:lap_i}. Moreover, consider again the case $\psi(u)=\frac{1}{2}u^2+ \gamma u$, $\gamma\geq0$, i.e.~$(X,\Q)$ is a
Bessel process of index $\gamma$.
\begin{enumerate}
\item Then, \eqref{eq:hw} corresponds to Hartman's result
\eqref{eq:hr}. Moreover, by letting $A \rightarrow \infty$ and using
the asymptotic behavior of the modified Bessel function of the first
kind, see \eqref{eq:bessel_asymp}, we get the Laplace
transform of the so-called Hartman-Watson law
\cite{Hartman-Watson-74}, the density of which has been
characterized by Yor \cite{Yor-80}.  Note also that this law is the mixture distribution in the representation of the Von-Mises distribution as a mixture of wrapped normal distributions, see \cite{Hartman-Watson-74}.

\item Moreover, by choosing $\gamma=\frac{1}{2}$, the mapping on the right hand side of \eqref{eq:id} corresponds to the Laplace transform of the L\'evy's
stochastic area integral. Indeed, for $B_t = (B^1_t ,B^2_t)$, a Brownian
motion on $\mathfrak{R}^2$, L\'evy \cite{Levy-50} computed the
Laplace transform of the process $L_t = \int_0^t B^1_s dB^2_s -B^2_s
dB^1_s, t>0,$ for fixed $u>0$, and $a = (\sqrt{u},\sqrt{u}) \in
\mathfrak{R}^2$, as follows
\begin{equation*}
 \E[e^{itL_u} |B_u = a] = \frac{tu}{\sinh(tu)} \exp{-(tu
\coth(tu) -1)}, \quad t \in \R.
\end{equation*}
This has been generalized by Biane and Yor
\cite{Biane-Yor-87a} to the L\'evy's stochastic area integral associated to some planar Gaussian Markov processes, in terms of the modified Bessel functions
of the first kind  of any index $\gamma > 0$.
\end{enumerate}
\end{remark}

\section{Proofs}
\subsection{Proof of Theorem \ref{thm:1}}
First, since the mapping $x\mapsto \Ip(x^{\alpha})$ is analytic on the right-half plane, it is plain that $\Ip \in  \mathbb{D}(\mathbf{L})$.  Observe also that, for any $\beta >0$, $x^{\beta} \in  \mathbb{D}(\mathbf{L})$ and
\begin{eqnarray} \label{eq:inf_lamp}\mathbf{L}x^{\beta} &=&
x^{\beta-\alpha}\psi(\beta).
\end{eqnarray}
Then, for any positive integer $N$, using the linearity of the operator
$\mathbf{L}$ and \eqref{eq:inf_lamp}, we get, for any $x\geq0 $,
\begin{eqnarray*}
\mathbf{L}\sum_{n=0}^{N} \ann q^nx^{\alpha n}&=&\sum_{n=1}^{N}
\ann q^n\mathbf{L}x^{\alpha n} \\
& = & \sum_{n=1}^{N} \ann q^n \psi(\alpha n)x^{\alpha(n-1)}\\
& = &q \sum_{n=0}^{N-1} \ann q^n x^{\alpha n}.
\end{eqnarray*}
The series being analytic on the right-half plane, then the right-hand side  of the previous line convergences as $N\rightarrow \infty$. Hence, we get  by monotone convergence (the series has only positive terms) that
\begin{eqnarray*}
\mathbf{L}\Ip(qx^{\alpha})&=& q\Ip(qx^{\alpha}), \quad x \geq 0 .
\end{eqnarray*}
Moreover, recalling that for $\E[\xi_1] <0 $, $\theta<\alpha$, we derive, in this case, that
\begin{eqnarray*}
\lim_{x \rightarrow 0} x^{-\theta+1}\frac{\partial}{\partial x}\Ip(x^{\alpha}) &=& \lim_{x \rightarrow 0} \sum_{n=0}^{\infty} \alpha(n+1)a_{n+1}(\psi, \alpha) x^{\alpha (n+1) -\theta}\\ &=& 0.
\end{eqnarray*}
Hence, for $\E[\xi_1] <0 $, $\Ip$ satisfy the condition \eqref{eq:cond-b}. Thus, it is an eigenfunction for the recurrent extension $(X,\Q)$.
Next, applying Dynkin's formula \cite{Dynkin-69} with the bounded stopping time $t \wedge \kappa_a$, we get, for any $t\geq0$,
\begin{eqnarray*}
  {\rm{E}}_x \left[e^{-q (t \wedge \kappa_a)}\Ip(qX^{\alpha}_{t \wedge \kappa_a}) \right]&=&
\Ip(qx^{\alpha}).
\end{eqnarray*}
Since the mapping $a\mapsto \Ip(qa)$ is increasing on $[0,\infty)$, we obtain, by dominated convergence and by using the fact that the process has no positive jumps, that, for any $0\leq x\leq a$,
\begin{eqnarray*}
  {\rm{E}}_x \left[e^{-q \kappa_a } \right]&=&
\frac{\Ip(qx^{\alpha})}{\Ip(qa^{\alpha})}.
\end{eqnarray*}
The proof of \eqref{eq:lap_ssm} is completed.
Next, the Esscher transform \eqref{eq:esscher} combined with the
Doob optional stopping theorem yields
\begin{eqnarray*}
\E_x \left[e^{-\lambda
\tau_a-q\Sigma_{\tau_a}}\Id{\tau_a<+\infty}\right]&=&
e^{\rho(x-a)}\E^{(\rho)}_x
\left[e^{-q\Ex_{\tau_a}}\Id{\tau_a<+\infty} \right].
\end{eqnarray*}
The proof of  Theorem \ref{thm:1} is completed by invoking the obvious identity $(\kappa_{e^a},\Q_{e^x}) \stackrel{(d)}{=}(\Sigma_{\tau_a},\P_x)$ on $\left\{\kappa_{e^a}<\kappa_{0}\right\}$ and by using \eqref{eq:lap_ssm}.
\hfill  $\square$

\vspace{0.5cm}

We end up this part by providing an interesting absolute continuity
relationship between self-similar processes, which will be useful for
the sequel. We recall from the discussion above that for $\E[\xi_1]
<0$, the self-similar process associated to $(\xi,\P)$ has a a.s. finite
lifetime $\kappa_0$. Then, by using the Lamperti's mapping
\eqref{eq:lamp_transf} and applying the chain rule, it follows that
the Esscher transform \eqref{eq:esscher}, reads, for $\gamma, \delta
\geq 0$ and $x >0$, as follows
\begin{equation}\label{eq:absl}
d\mathbb{Q}^{(\gamma)}_{x|F_{A_t}} =
   \left(\frac{X_t}{x}\right)^{\gamma-\delta} e^{- \left(\psi(\gamma)-\psi(\delta)\right) \iex_t }   \,
   d\mathbb{Q}^{(\delta)}_{x|F_{A_t}},\quad \textrm{ on }
   \{t<\kappa_0\},
\end{equation}
where
\begin{eqnarray*}
 A_t =  \int_0^{t} X^{-\alpha}_u  \d u.
\end{eqnarray*}
Since for every $F_{A_t}$-stopping time $T$, $\Sigma(T)$ is a
$F_{t}$-stopping time,  the absolute continuity relationship
\eqref{eq:absl} holds for every $F_{A_t}$-stopping time on
$F_{A_{T^+}} \cup \{T<\kappa_0\}$.  We mention that, for
$\E[\xi_1]>0$, where the condition "on $\{t<\kappa_0\}$" can be
omitted, the relationship \eqref{eq:absl} was already established by
Carmona et al.~\cite[Proposition 2.1]{Carmona-Petit-Yor-94}, under
the name of Girsanov power transformation. Finally, if we set
$\gamma=\theta$ and $\delta=0$ in the case $\E[\xi_1] <0$, then
\eqref{eq:absl} simplifies to the following Doob's $h$-transform,
for $x>0$,
\begin{equation}\label{eq:absl0}
d\mathbb{Q}^{(\theta)}_{x|F_{A_t}} =
   \left(\frac{X_t}{x}\right)^{\theta}d\mathbb{Q}_{x|F_{A_t}},\quad \textrm{ on } \{t<\kappa_0\}.
\end{equation}

\subsection{Proof of Theorem \ref{thm:2}}
We assume that $\E[\xi_1] <0$ and $\theta<\alpha$. Then,  the identity  $(x^{\alpha}\Ex_{\infty}, \P)
\stackrel{(d)}{=}(\kappa_{0}, \Q_x)$ yields
\begin{eqnarray*}
{\rm{E}}_x \left[e^{-q \kappa_0}\right]&=& {\E}\left[e^{-qx^{\alpha} \Sigma_{\infty}}\right].
\end{eqnarray*}
Hence, the mapping $x \mapsto {\rm{E}}_x \left[e^{-q \kappa_0}\right]$ is decreasing on $[0,\infty)$,  and  by dominated convergence we have, see also Vuolle-Apiala \cite{Vuolle-Apiala-94},
\begin{eqnarray}
\lim_{x \rightarrow
\infty}{\rm{E}}_x \left[e^{-q \kappa_0}\right]&=& \lim_{x \rightarrow
\infty}{\E}\left[e^{-qx^{\alpha} \Sigma_{\infty}}\right] \nonumber \\
&=& 0  \label{eq:lim0}
\end{eqnarray}
and
\begin{eqnarray*}
\lim_{x \rightarrow
0}{\rm{E}}_x \left[e^{-q \kappa_0}\right]&=& 1.
\end{eqnarray*}
We now compute the Laplace transform of $\kappa_0$. From the
Doob's $h$-transform \eqref{eq:absl0}, we deduce that, for any $x,a$ such
that $0<x\leq a$,
\begin{eqnarray*}
{\rm{E}}_x \left[e^{-q \kappa_{a}}\Id{\kappa_a<\kappa_{0}} \right]
&=&\left(\frac{x}{a}\right)^{\theta}\frac{\It\left(qx^{\alpha}\right)}{\It\left(qa^{\alpha}\right)}.
\end{eqnarray*}
Then, the strong Markov property and the absence of positive jumps yield
\begin{eqnarray*}
{\rm{E}}_x \left[e^{-q \kappa_{0}}\Id{\kappa_{0}<\kappa_a} \right]
&=&\frac{1}{{\rm{E}}_0 \left[e^{-q \kappa_a}\right]}\left( {\rm{E}}_x \left[e^{-q \kappa_a}\right]-{\rm{E}}_x \left[e^{-q \kappa_a}\Id{\kappa_a<\kappa_{0}} \right]\right)\\
&=&\Ip\left(qx^{\alpha}\right)-\Ip\left(qa^{\alpha}\right)
\left(\frac{x}{a}\right)^{\theta}\frac{\It\left(qx^{\alpha}\right)}{\It\left(qa^{\alpha}\right)}
\end{eqnarray*}
and
\begin{eqnarray*}
{\rm{E}}_x \left[e^{-q \kappa_{0}} \right] &=&{\rm{E}}_x \left[e^{-q
\kappa_{0}}\Id{\kappa_{0}<\kappa_a} \right]+{\rm{E}}_x
\left[e^{-q \kappa_a} \Id{\kappa_a<\kappa_{0}}\right]{\rm{E}}_a \left[e^{-q \kappa_{0}} \right]\\
&=& \Ip\left(qx^{\alpha}\right)-\Ip\left(qa^{\alpha}\right)
\left(\frac{x}{a}\right)^{\theta}\frac{\It\left(qx^{\alpha}\right)}{\It\left(qa^{\alpha}\right)}\\&+&\left(\frac{x}{a}\right)^{\theta}\frac{\It\left(qx^{\alpha}\right)}{\It\left(qa^{\alpha}\right)}{\rm{E}}_a
\left[e^{-q \kappa_{0}} \right]\\
&=& \Ip\left(qx^{\alpha}\right)-x^{\theta}
\It\left(qx^{\alpha}\right)
\frac{a^{-\theta}}{\It\left(qa^{\alpha}\right)}\left(
\Ip\left(qa^{\alpha}\right)-{\rm{E}}_a
\left[e^{-q \kappa_{0}}
\right]\right).
\end{eqnarray*}
To derive an expression of the sought quantity, we first
differentiate with respect to $a$ the previous equation, the
function involved being smooth, to get the Riccati equation
\begin{eqnarray*}
\frac{a^{-\theta}}{\It\left(qa^{\alpha}\right)} \left(\frac{\partial}{\partial a}{\rm{E}}_a
\left[e^{-q \kappa_{0}}\right] -\frac{\partial}{\partial a}\Ip\left(qa^{\alpha}\right)\right) + \left({\rm{E}}_a
\left[e^{-q \kappa_{0}}\right] - \Ip\left(qa^{\alpha}\right)\right)\frac{\partial}{\partial a}
\frac{a^{-\theta}}{\It\left(qa^{\alpha}\right)} &=&0
\end{eqnarray*}
where we have used the fact that $x^{\theta}
\It\left(qx^{\alpha}\right)>0$ for any $x>0$. By
solving this Riccati equation, we observe that the
solution has the following form
\begin{eqnarray} \label{eq:h0}
{\rm{E}}_x \left[e^{-q \kappa_{0}} \right]
&=&A\Ip\left(qx^{\alpha}\right)-Cx^{\theta}\It\left(qx^{\alpha}\right)
\end{eqnarray}
for some constants $A,C$. It is immediate that $A=1$ since
${\rm{E}}_0 \left[e^{-q \kappa_0}\right]=1$. Then, the
self-similarity property yields $C=q^{\frac{1}{\alpha}}c$, for some
real constant c. Furthermore, \eqref{eq:lim0}
 ensures the
existence of a constant $C_{\theta}>0$ such that
\begin{equation*}
\Ip\left(x^{\alpha}\right)\sim C_{\theta}
 x^{\theta}\It\left(x^{\alpha}\right) \quad {\rm{as }} \: x\rightarrow \infty.
 \end{equation*}
Moreover, observe, from \eqref{eq:h0}, that
\begin{eqnarray*}
\lim_{x\rightarrow 0}\frac{{\rm{E}}_x \left[1-e^{-\kappa_{0}}
\right]}{x^{\theta}} = C_{\theta}.
\end{eqnarray*}
Next,  we recall the following identities, see \cite[(11) and Remarks 1. p.~489]{Rivero-05},
\begin{eqnarray*}
\lim_{x\rightarrow 0}\frac{{\rm{E}}_x \left[1-e^{-\kappa_0}
\right]}{x^{\theta}} &=&
\frac{\Gamma(1-\frac{\theta}{\alpha})}{\alpha
\psi'(\theta)}\E\left[\Sigma_{\infty}^{\frac{\theta}{\alpha}-1}\right]
\end{eqnarray*}
and
\begin{eqnarray*}
\E\left[\Sigma_{\infty}^{\frac{\theta}{\alpha}-1}\right]&=&
\E^{(\theta)}\left[\left(\int_0^{\infty}e^{-\alpha
\xi_s}ds\right)^{\frac{\theta}{\alpha}-1}\right].
\end{eqnarray*}
Since $(\xi,\P^{(\theta)})$ has a positive mean,  the proof of the
 Theorem \ref{thm:2} is completed by using Proposition 2.3 in
Maulik and Zwart \cite{Maulik-Zwart-06} and by choosing $c =
C_{\theta}$. \hfill  $\square$

\subsection{Proof of Corollary \ref{cor:1}} We characterize the entrance laws of $(\widehat{X},\Q)$ and
$(\widehat{X},\Q^{(\theta)})$. To this end, we state the following easy result.
\begin{lemma}
Let us assume that $-\infty <\E[\xi_1] <0$. Then, the entrance law of $(\widehat{X},\Q)$ admits a density $\eta_t(y)$  with respect to the reference measure $m(dy), \: y>0$. Moreover, the Laplace transform in time of $\eta_t(y)$, denoted by $n^{q}(y),\: q \geq 0$, is characterized by the identity
\begin{equation*}
n^q(y) = \frac{1}{|\E[\xi_1]|} \E\left[e^{-qy^{\alpha}\Sigma_{\infty}}\right].
\end{equation*}
\end{lemma}
\begin{proof}
The claims follow readily from Bertoin and Yor \cite[p.~396]{Bertoin-Yor-02}. Indeed, they characterize the $q$-potential of the
entrance law of a self-similar process associated, via the Lamperti's
mapping, to a L\'evy process with a positive and finite mean,  as
follows. For any measurable function $f : \R^+\rightarrow \R^+$, we have
\begin{equation*}
n^q f = \frac{1}{|\E[\xi_1]|}\int_0^{\infty} f(y) \E\left[e^{-qy^{\alpha}\Sigma_{\infty}}\right] m(dy)
\end{equation*}
where we have used the identity $\E[\hat{\xi}_1]=-\E[\xi_1]$.
\end{proof}
The first part of the Corollary follows from \eqref{eq:lt_exp}.
For the second statement, we use \cite[Proposition 3]{Rivero-05} where the $q$-potential, $n^q$, of the entrance law is given, for a bounded continuous function $f$, by
\begin{equation*}
n^q f = \frac{1}{\psi'(\theta) C_{\theta}}\int_0^{\infty} f(y) \E\left[e^{-qy^{\alpha}\Sigma_{\infty}}\right] y^{\alpha-1-\theta}dy
\end{equation*}
The proof of the Corollary is completed by means of the identity \eqref{eq:lt_exp}.

Note the following result regarding the resolvent of the minimal process $(X,\kappa_0)$.
\begin{lemma}
For $\E[\xi_1] <0$, the resolvent, $U_0^{q}$, of the minimal process $(X,\kappa_0)$ admits a density with respect to the reference measure $m(dy)$, which is jointly continuous and bounded. Thus, the processes $(X,\kappa_0)$ and $(\widehat{X},\Q)$ are in classical duality.
\end{lemma}
\begin{proof}
We recall that in the case $\E[\xi_1] <0$, the process $(\xi,\P)$ is necessarily of unbounded variation since we have excluded the case of negative subordinators. Thus, each point of the real line is regular for itself, see  \cite[VII, Corollary 5]{Bertoin-96}.  Moreover, since the $q$-capacity of $\{0\}$, which is $\phi'(q)$, is positive for any $q>0$, see \cite[VII.5.2]{Bertoin-96}, and the resolvent of $(\xi,\P)$ is absolutely continuous, we deduce that points are not polar for $(\xi,\P)$, i.e.~$\P_x(\tau_y<+\infty)>0$ for any $x,y \in \R$. It is not difficult
to see that these two properties are left invariant by time change
with a continuous additive functional, see Bally and Stoica \cite[Proposition 4.1]{Bally-Stoica-87}. The assertions follow from \cite[Proposition 3.1]{Bally-Stoica-87}.
\end{proof}

\subsection{Proof of Theorem \ref{thm:3}}
The claim
\eqref{eq:selfd} is contained in  the following Lemma.
\begin{lemma} \label{lem:sd}
The process $(\kappa_a)_{a \geq 0}$ is under $\Q_{0^+}$ an
$\alpha$-self-similar additive process, i.e.~a
 process with independent increments which enjoy the scaling property \eqref{eq:self}.
Hence $\kappa_1$ is under $\Q_{0^+}$ a positive self-decomposable random variable.
\end{lemma}
\begin{proof}
The first assertion follows from the absence of positive jumps, the strong Markov property and the
self-similarity of $(X, \Q)$.
 The last statement is a
straightforward consequence of the property of the law of additive
processes,  see Sato \cite[Chap. 3.16]{Sato-99} .
\end{proof}
Moreover, it is well know, see Wolfe \cite{Wolfe-82} and also
Jeanblanc et al.~\cite{Jeanblanc-Pitman-Yor-02} for related results, that  if the mapping $q \mapsto f(q)$ is the Laplace transform of a positive self-decomposable random variable then there exists a unique, in distribution, (increasing) L\'evy process
$L$ such that $\E\left[\log(1+L_1)\right] < +\infty$ and its Laplace
exponent, $\phi_L$, is given by
\begin{equation*}
\phi_L(q) = q\frac{\frac{d}{dq}f(q)}{f(q)}, \: q\geq0.
\end{equation*}
Thus, from Theorem \ref{thm:2} and Lemma \ref{lem:sd}, we deduce the statement \eqref{eq:id},  after recalling that the infinite divisibility property is stable under convolution.

 Before stating the second
Lemma, we assume $\E[\xi_1]\geq 0$ and we introduce some notation. Let $(P_t)_{t\geq0}$ be the semigroup of the L\'evy process $\xi$. We denote by $(P^{\Sigma}_t)_{t\geq0}$ the subordinate semigroup of $(P_t)_{t\geq0}$ by the continuous decreasing multiplicative functional $(e^{-q\Sigma_t})_{t\geq0}$.  That is for $f \in B(\R)$, we have, for any $t\geq0$,
\begin{equation*}
P^{\Sigma}_t f(x) = \E_x\left[e^{-q\Sigma_t} f(\xi_t)\right], \: x \in \R.
\end{equation*}
 Next, by choosing $\lambda=0$ in \eqref{eq:jlp}, we deduce that the function  $x \mapsto \Ip\left(qe^{\alpha x}\right)$ is excessive for the semigroup
$(P^{\Sigma}_t)_{t\geq0}$. Moreover, it is plain that, for any $x\in \R$, $0<\Ip\left(qe^{\alpha x}\right)<\infty$. Thus, one can define a new real-valued (sub)-Markov process with semigroup (resp.~ law)
denoted by $(P_t^{\I})_{t\geq0}$ (resp.~$\P^{\I}$), as Doob's $h$-transform of $(P^{\Sigma}_t)_{t\geq0}$, as follows, for any $f \in B(\R)$ and $t\geq 0$,
\begin{equation} \label{eq:eq_h}
P^{\I}_t f(x) = \frac{1}{\Ip\left(qe^{\alpha x}\right)}P^{\Sigma}_t \left(f\Ip\left(qe^{\alpha.}\right)\right)(x), \: x \in \R.
\end{equation}
 We are now
ready to state the following result which characterizes the random
variable associated to the Laplace transform \eqref{eq:hw}.
\begin{lemma}
For any $0\leq x \leq a $, $\lambda \geq 0$ and recalling that $\rho=\phi(\lambda)$ , we have
\begin{equation} \label{eq:lap_i}
\E_x^{\I}\left[e^{-\lambda \tau_a }\right]
=e^{\rho(x-a)}\frac{\Ia\left(qe^{\alpha
x}\right)\Ip\left(qe^{\alpha a}\right)}{\Ia\left(qe^{\alpha a}\right)\Ip\left(qe^{\alpha x}\right)}.
\end{equation}
\begin{proof}
We deduce from \eqref{eq:eq_h} that the following absolute continuity
relationship
\begin{equation}\label{eq:eq_ph}
d\mathbb{P}^{\I}_{x|F_{t}} =
   \frac{e^{-q\Sigma_t}\Ip\left(qe^{\alpha \xi_t}\right)}{\Ip\left(qe^{\alpha x} \right)}   \,
   d\mathbb{P}_{x|F_{t}}
\end{equation}
holds for any $t>0$ and $x\in\R$. It is plain that this relationship remains
valid on $F_{{T^+}} \cap \{T^+ < \infty\}$ for any
$F_{\infty}$-stopping time $T$.  Recalling that  $\P_x\left[\tau_a<\infty\right]=1$ for $\E[\xi_1] \geq 0$, we get that
\begin{eqnarray*}
\E_x^{\I}\left[e^{- \lambda \tau_a } \right]  &=&
\frac{\Ip\left(qe^{\alpha a}\right)}{\Ip\left(qe^{\alpha x}\right)}\E_x \left[e^{- \lambda \tau_a - q\Sigma_{\tau_a}} \right] \\
&=& e^{\rho(x-a)}\frac{\Ia\left(qe^{\alpha
x}\right)\Ip\left(qe^{\alpha a}\right)}{\Ia\left(qe^{\alpha a}\right)\Ip\left(qe^{\alpha x}\right)}
\end{eqnarray*}
where the last line follows from \eqref{eq:jlp}.
\end{proof}
\end{lemma}
Finally, it is plain from the absolute continuity \eqref{eq:eq_ph} that the process under $\P^{\I}$ is also spectrally negative. Then, from the strong Markov and the absence of positive jumps, we have for any $x< c < a$,
\begin{equation*}
    \left(\tau_{a},\P_x^{\I}\right) \stackrel{(d)}{=}  \left(\tau_{c},\P_x^{\I}\right) + \left(\tau_{a},\P_c^{\I}\right)
\end{equation*}
where the random variable on the right-hand side are independent.  Hence, $\left(\tau_{a},\P_x^{\I}\right)$ is infinitely divisible.
 The proof of Theorem \ref{thm:3} is then completed.\hfill  $\square$

\section{Some illustrative examples}
We end up by investigating some well-known and new examples in more
details.

\subsection{The modified Bessel functions} \label{sec:b}
We consider $\xi$ to be a Brownian motion with drift $\gamma \in \R$, i.e.~$\psi(u)=\frac{1}{2}u^2+\gamma u$ and we set $\alpha=2$. In the case $\gamma<0$, we have
$\theta=2\gamma$ and therefore we assume $\gamma>-1$. Its associated
self-similar process is well known to be a Bessel process of index
$\gamma$. In the sequel, we simply indicate the connections between the
power series $\I_{2,\psi}$ and the modified Bessel functions since the
results of this paper are well known and can be found for instance
in Hartman \cite{Hartman-76}, Pitman and Yor \cite{Pitman-Yor-80}
and Yor \cite{Yor-01}. We have
\begin{eqnarray*} a_n(\psi,\gamma;2)^{-1}&=&2^{n} n!\prod_{k=1}^n
\left(
k+\gamma\right)\\
&=& 2^{n} n! \frac{\Gamma\left( n+ \gamma +1\right)}{\Gamma(\gamma+1)},
\quad a_0=1.
\end{eqnarray*}
Thus, we get
\begin{equation*}
\I_{2,\psi}(x)=(x/2)^{-\gamma/2}\Gamma(\gamma+1){\rm{I}}_{\gamma}\left(\sqrt{2x}\right)
\end{equation*}
where
\begin{equation*}
{\rm{I}}_{\gamma}(x)=\sum_{n=0}^{\infty}\frac{(x/2)^{\gamma+2n}}{n!\Gamma(\gamma+n+1)}
\end{equation*}
 stands for the modified Bessel function of index $\gamma$, see e.g.~\cite[Chap.~5]{Lebedev-72}. The asymptotic behavior of this function is well known to
\begin{equation} \label{eq:bessel_asymp}
{\rm{I}}_{\gamma}(x) \sim \frac{e^{x}}{\sqrt{2\pi x}} \quad  \textrm{ as } x
\rightarrow \infty.
\end{equation}
Thus, we obtain that
 \begin{eqnarray*}
\N_{2,\psi_{2\gamma}}(x)&=&(x/2)^{\gamma/2}\Gamma(-\gamma+1)\left({\rm{I}}_{-\gamma}\left(\sqrt{2 x}\right)-{\rm{I}}_{\gamma}\left(\sqrt{2 x}\right)\right)\\
&=&(x/2)^{\gamma/2}\frac{2}{\Gamma(\gamma)}{\rm{K}}_{\gamma}\left(\sqrt{2x}\right)
\end{eqnarray*}
where $2 {\rm{K}}_{\gamma}(x) = \Gamma(1-\gamma) \Gamma(\gamma)
({\rm{I}}_{-\gamma}(x)-{\rm{I}}_{\gamma}(x))$ is the MacDonald function of
index $\gamma$.
\subsection{Some generalizations of the Mittag-Leffler function} \label{sec:ml}
In \cite{Patie-06-poch}, the author introduced a new parametric
family of one-sided L\'evy processes which are characterized by the
following Laplace exponent, for any  $ 1<\varrho< 2$, $\beta \geq0$
and $\gamma
>1-\varrho$,
\begin{equation} \label{eq:lap_poch}
 \psi(\beta u+\gamma)-\psi(\gamma)=
\frac{1}{\varrho}\left((\beta
u+\gamma-1)_{\varrho}-(\gamma-1)_{\varrho}\right)
\end{equation}
where $(k)_{\varrho}=\frac{\Gamma(k+\varrho)}{\Gamma(k)}$ stands for
the Pochhammer symbol. Its  characteristic triplet  is $\sigma=0$,
\begin{equation*} \tilde{\nu}(dy)=\frac{ \varrho(\varrho-1)}{\beta\Gamma(2-\varrho)
}\frac{e^{(\varrho+\gamma-1)\frac{y}{\beta}}}{(1-e^{\frac{y}{\beta}})^{\varrho+1}}dy,
\quad y<0,
\end{equation*}
and
\begin{equation*} b_{\gamma}=\beta
(\gamma)_{\varrho}(\Psi(\gamma-1+\varrho)-\Psi(\gamma-1))
\end{equation*}
where $\Psi(\lambda)=\frac{\Gamma'(\lambda)}{\Gamma(\lambda)}$ is
the digamma function. In particular, if $\gamma_0 $ denotes  the
zero of the function $\gamma \rightarrow b_{\gamma}  $, then for
$\gamma \geq\gamma_0\in (1-\varrho, 0)$, $\E[\xi_1]\geq 0$.
 \subsubsection{The case $\gamma=0$}\eqref{eq:lap_poch} reduces to $ \psi(u)=
\frac{1}{\varrho}(u-1)_{\varrho}$. Observe that $\theta=1$ and $\psi'(1)=\frac{\Gamma(\varrho)}{\varrho}$. Moreover, setting $\alpha=\varrho$, we get
\begin{equation*}
a_n(\varrho,0;\varrho)^{-1}=\frac{\Gamma(\varrho(n+1)-1)}{\Gamma(\varrho-1)},
\quad a_0=1.
\end{equation*}
The series, in this case, can be written as follows
\begin{eqnarray*}
\I_{\varrho,\psi_{\theta}}(x)&=&\Gamma(\varrho)\M_{\varrho,\varrho}(\varrho
x)\\
\I_{\varrho,\psi}(x)&=&\Gamma(\varrho-1)\M_{\varrho,\varrho-1}(\varrho
x)
\end{eqnarray*}
where
\begin{equation*}
 \M_{\varrho,\beta}(z)= \sum_{n=0}^{\infty}  \frac{ z^n}{\Gamma(\varrho n+\beta)}, \quad
 z
\in  \C,
\end{equation*}
stands for the Mittag-Leffler function of parameter
$\varrho,\beta>0$. The function $\M_{\varrho,0}(z)$ was defined and
studied by Mittag-Leffler \cite{Mittag-03}. It is a direct
generalization of the exponential function. The generalization
$\M_{\varrho,\beta}(z)$ was given by Agarwal \cite{Agarwal-53}
following the work of Humbert \cite{Humbert-53}. A detailed account
of these functions is available from the monograph of Erd\'elyi et
al.~\cite{Erdelyi-55}. Next, we recall the following Mellin-Laplace
transform of the generalized Mittag-Leffler function, for
$\lambda,\beta \geq0$, we have
\begin{eqnarray*}
\int_{0}^{\infty}e^{-(\lambda+1)x}x^{\beta-1}\M_{\varrho,\beta}(x^{\varrho})dx
&=& \frac{(\lambda+1)^{\varrho-\beta}}{(\lambda+1)^{\varrho}
-1}.\end{eqnarray*}
 We deduce  by invoking a
Tauberian theorem, see Feller \cite[XIII.5]{Feller-71}, the
following asymptotic behavior
\begin{equation*}
 \M_{\varrho,\beta}(x^{\varrho})\sim
\frac{1}{\varrho}e^{x}x^{1-\beta} l(x^{\varrho}) \quad {\rm{ as }} \:
x \rightarrow \infty,
\end{equation*}
with $l$ a slowly varying function at infinity. Thus,
$C_{1}=\frac{\varrho}{\varrho-1}$,
$\E\left[\Sigma_{\infty}^{\frac{1}{\varrho}-1}\right] = \frac{\varrho
\Gamma(\varrho-1)}{\Gamma(1-\frac{1}{\varrho})}$ and
\begin{equation*}
\N_{\varrho,\psi_1}(x^{\varrho}) =
\M_{\varrho,\varrho-1}(x^{\varrho})-\frac{\varrho
x}{\varrho-1}\M_{\varrho,\varrho}(x^{\varrho}).
\end{equation*}
Finally, we can state  the following properties of the
Mittag-Leffler function.
\begin{cor}
Let $1<\varrho< 2$. The mappings
\begin{eqnarray*}
 q \mapsto \frac{1}{\M_{\varrho,\varrho}(q)} \quad {\textrm{ and }} \quad q \mapsto  \M_{\varrho,\varrho-1}(q)-\frac{\varrho q^{\frac{1}{\varrho}}}{\varrho-1}\M_{\varrho,\varrho}(q)
\end{eqnarray*}
are respectively the Laplace transform of a positive self-decomposable
distribution and a completely monotone function. The mapping
\begin{eqnarray*}
 q \mapsto
 \exp\left(\frac{q\M'_{\varrho,\varrho}(q)}{\M_{\varrho,\varrho}(q)}\right)
\end{eqnarray*}
is the Laplace transform of a positive infinitely divisible
distribution.
\end{cor}

\subsubsection{The general case} In this case,
\begin{equation*}
a_n(\beta,\gamma;\varrho)^{-1}=\prod_{k=1}^n \left(( \beta k
+\gamma-1)_{\varrho}-(\gamma-1)_\varrho\right), \quad a_0=1.
\end{equation*}
and write $\M_{\varrho,\beta,\gamma}(x) =
\sum_{n=0}^{\infty}a_n(\beta,\gamma;\varrho)x^{n}$. We point out that
this function is closely related to the power series introduced by
Kilbas and Saigo \cite{Kilbas-Saigo-95} which have coefficients of
the following form
\begin{equation*}
\tilde{a}_n( \beta,\gamma;\varrho)^{-1}=\prod_{k=1}^n (\varrho( \beta
k +\gamma)+1)_{\varrho}, \quad a_0=1.
\end{equation*}
From Theorem \ref{thm:3}, we deduce the following properties.
\begin{cor}
Let  $1<\varrho < 2$, $\beta \geq 0$ and $\gamma>1-\varrho$. Then,
the mapping
\begin{eqnarray*}
 q \mapsto \frac{1}{\M_{\varrho, \varrho \beta, \gamma}(q)}
\end{eqnarray*}
is the Laplace transform of a positive self-decomposable
distribution. The mapping
\begin{eqnarray*}
 q \mapsto \exp\left( \frac{q\M'_{\varrho, \varrho \beta, \gamma}(q)}{\M_{\varrho, \varrho \beta,
 \gamma}(q)}\right)
\end{eqnarray*}
is the Laplace transform of a positive infinitely divisible
distribution. Finally, for $0<a< A$ and $\gamma\geq \gamma_0$, recalling that $\rho=\phi(\lambda)$, the
mapping
\begin{eqnarray*}
  \lambda \mapsto \left(\frac{a}{A}\right)^{\rho} \frac{\M_{\varrho, \varrho \beta, \gamma+\rho}(a)\M_{\varrho, \varrho \beta, \gamma}(A)}{\M_{\varrho, \varrho \beta, \gamma+\rho}(A)\M_{\varrho, \varrho \beta, \gamma}(a)}
\end{eqnarray*}
is the Laplace transform of a positive infinitely divisible
distribution.
\end{cor}
\subsection{The power series associated to stable processes and a new generalization of the exponential function}
Finally, we consider the Esscher transform of spectrally negative
stable process, i.e. $\psi_{\gamma}(u) =
c_{\varrho}\left((u+\gamma)^{\varrho}-(\gamma)^{\varrho}\right)$,
$\gamma\geq0$, $1<\varrho <2$ and  $c_{\varrho}>0$. Its
characteristic triplet is $\sigma=0,$
\begin{equation*} \tilde{\nu}(dy)=\frac{ c_{\varrho}\varrho(\varrho-1)}{\Gamma(2-\varrho)
}\frac{e^{\gamma y}}{\mid y\mid^{\varrho+1}}dy, \quad y<0,
\end{equation*}
and $b=c_{\varrho}\varrho\gamma^{\varrho-1} \geq0$. The inverse
function of $\psi$ is $\phi(u)=(c_{\varrho}
u+\gamma^{\varrho})^{\frac{1}{\varrho}}-\gamma$ and
\begin{equation*}
a_n(\alpha,\varrho,\gamma)^{-1}= \prod_{k=1}^n \left(
\left(\alpha k+\gamma \right)^{\varrho}-\gamma^{\varrho}\right), \quad
a_0=1.
\end{equation*}
Such formulation motivates us to introduce a generalization of the
factorial symbol,
 which we defined, for $n \in \mathfrak{N}$ ,$\alpha\in \C,
\Re(\alpha)\geq0$  and $\gamma\in \C, \Re(\gamma)\geq0$, by
\begin{equation*}
 (\alpha,\gamma)_{\varrho,n} = \prod_{k=1}^n \left((k\alpha+\gamma)^{\varrho}
 -\gamma^{\varrho}\right)
\: \textrm{ and } (\alpha,\gamma)_{\varrho,0} =1.
\end{equation*}
Note the obvious identities
\begin{eqnarray*}
(\alpha,\gamma)_{0,n}  &=& 0 \\
(\alpha,\gamma)_{1,n}  &=& \alpha^{n}n! \\
(\alpha,0)_{\varrho,n}  &=& \alpha^{n} n!(\alpha,0)_{\varrho-1,n} \\
(\alpha,\gamma)_{\varrho,n} &=& \alpha^{\varrho
n}\left(1,\frac{\gamma}{\alpha}\right)_{\varrho,n} .
\end{eqnarray*}
Moreover, we have
\begin{eqnarray*}
(\alpha,\gamma)_{\varrho,n}&=& \alpha^{\varrho n}\sum_{k=0}^n (-1)^{n-k}
\left(\frac{\gamma}{\alpha}\right)^{\varrho(n-k)}\left(1,\frac{\gamma}{\alpha}\right)_{\varrho,k}.
\end{eqnarray*}
We write simply
\begin{equation*}
\I_{\alpha,\gamma,\varrho}(c_{\varrho}x)=\sum_{n=0}^{\infty}
\frac{x^n}{(\alpha,\gamma)_{\varrho,n}}.
\end{equation*}
Observe that
\begin{equation*}
\lim_{\varrho \downarrow 1}
\I_{\alpha,\gamma,\varrho}(c_{\varrho}x)=e^{\frac{x}{\alpha}}.
\end{equation*}
\begin{cor}
Let $1<\varrho < 2$, $\alpha>0$ and $\gamma \geq 0$. Then, the mapping
\begin{eqnarray*}
 q \mapsto \frac{1}{\I_{\alpha,\gamma,\varrho}(q)}
\end{eqnarray*}
is the Laplace transform of a positive self-decomposable
distribution. The mapping
\begin{eqnarray*}
 q \mapsto \exp\left(
 \frac{q\I'_{\alpha,\gamma,\varrho}(q)}{\I_{\alpha,\gamma,\varrho}(q)}\right)
\end{eqnarray*}
is the Laplace transform of a positive infinitely divisible
distribution. Finally, for $0<a< A$ and writing $\rho=\phi(\lambda)$, the mapping
\begin{eqnarray*}
  \lambda \mapsto \left(\frac{a}{A}\right)^{\rho} \frac{\I_{\alpha,\rho,\varrho}(a)\I_{\alpha,0,\varrho}(A)}{\I_{\alpha,\rho,\varrho}(A)\I_{\alpha,0,\varrho}(a)}
\end{eqnarray*}
is the Laplace transform of a positive infinitely divisible
distribution.
\end{cor}


\begin{thebibliography}{10}

\bibitem{Agarwal-53}
R.~P.~Agarwal, \emph{A propos {d'une} note de {M. Pierre Humbert}}, C. R. Math.
  Acad. Sci. Paris \textbf{236} (1953), 2031--2032.

\bibitem{Bally-Stoica-87}
V.~Bally and L.~Stoica, \emph{A class of {Markov} processes which admit a local
  time}, Ann. Probab. \textbf{15} (1987), no.~1, 241--262.

\bibitem{Bertoin-96}
J.~Bertoin, \emph{L\'evy processes}, Cambridge University Press, Cambridge,
  1996.

\bibitem{Bertoin-Yor-02-b}
J.~Bertoin and M.~Yor, \emph{The entrance laws of self-similar {M}arkov
  processes and exponential functionals of {L}\'evy processes}, Potential Anal.
  \textbf{17} (2002), no.~4, 389--400.

\bibitem{Bertoin-Yor-02}
\bysame, \emph{On the entire moments of self-similar {M}arkov processes and
  exponential functionals of { L}\'evy processes}, Ann. Fac. Sci. Toulouse
  Math. \textbf{11} (2002), no.~1, 19--32.

\bibitem{Bertoin-Yor-05}
\bysame, \emph{{Exponential functionals of L\'evy processes}}, Probab. Surv.
  \textbf{2} (2005), 191--212.

\bibitem{Biane-Yor-87a}
Ph.~Biane and M.~Yor, \emph{Variations sur une formule de {Paul L\'evy}}, Ann.
  Inst. H. Poincar\'e Probab. Statist. \textbf{23} (1987), no.~2, 359--377.

\bibitem{Blumenthal-83}
R.M.~Blumenthal,  \emph{On construction of {M}arkov processes}, Z. Wahrsch. Verw. Gebiete. \textbf{63} (1983), 433--444.


\bibitem{Carmona-Petit-Yor-94}
Ph.~Carmona, F.~Petit, and M.~Yor, \emph{Sur les fonctionnelles exponentielles
  de certains processus de {L\'e}vy}, Stoch. Stoch. Rep. \textbf{47} (1994),
  no.~1-2, 71--101, (english version in \cite{Yor-01} p. 139--171).

\bibitem{Caballero-Chaumont-06}
M.E.~Caballero, and L.~Chaumont,
\emph{Weak convergence of positive self-similar {M}arkov processes
              and overshoots of {L}\'evy processes}, Ann. Probab.,
\textbf{34} (2006), no.~3, 1012--1034.

\bibitem{Ciesielski-Taylor-62}
Z.~Ciesielski and S.J.~Taylor.
\emph{First passage times and sojourn times for {B}rownian motion in space
  and the exact {H}ausdorff measure of the sample path}.
 Trans. Amer. Math. Soc. \textbf{103} (1962), 434--450.

\bibitem{Dynkin-69}
E.B.~Dynkin,  \emph{Markov processes. {V}ols. {I}, {II}}, vol.~122, Academic Press Inc., Publishers, New York,
  1969.
\bibitem{Erdelyi-55}
A.~Erd{\'e}lyi, W.~Magnus, F.~Oberhettinger, and F.G.~Tricomi, \emph{Higher
  transcendental functions}, vol.~3, McGraw-Hill, New York-Toronto-London,
  1955.

\bibitem{Feller-71}
W.E.~Feller, \emph{An introduction to probability theory and its applications},
  $2^{nd}$ ed., vol.~2, Wiley, New York, 1971.


\bibitem{Gikhman-Skorohod-75}
I.I.~Gikhman  and A.V.~Skorokhod, \emph{The theory of stochastic processes},
  springer- verlag, new york-heidelberg-berlin ed., vol.~II, 1975.


\bibitem{Hartman-76}
P.~Hartman, \emph{Completely monotone families of solutions of $n$-th order
  linear differential equations and infinitely divisible distributions}, Ann.
  Sc. Norm. Super. Pisa Cl. Sci. \textbf{IV-III} (1976), 267--287.

\bibitem{Hartman-Watson-74}
P.~Hartman and G.S.~Watson, \emph{"{N}ormal" distribution functions on spheres
  and the modified {B}essel functions}, Ann. Prob. \textbf{2} (1974), 593--607.

\bibitem{Humbert-53}
P.~Humbert, \emph{Quelques r\'esultats relatifs \`a la fonction de
  {Mittag-Leffler}}, C. R. Math. Acad. Sci. Paris \textbf{236} (1953),
  1467--1468.

\bibitem{Jeanblanc-Pitman-Yor-02}
M.~Jeanblanc, J.~Pitman and M.~Yor, \emph{{Self-}similar processes with
  independent increments associated with {L}\'evy and {Bessel} processes.},
  Stochastic Process. Appl. \textbf{100} (2002), 223--232.

\bibitem{Kent-78}
J.~Kent, \emph{Some probabilistic properties of {B}essel functions}, Ann.
  Probab. \textbf{6} (1978), no.~5, 760--770.

\bibitem{Kilbas-Trujilo-01}
A.A.~Kilbas and J.J.~Trujillo, \emph{Differential equations of fractional
  orders: Methods, results and problems}, Appl. Anal. \textbf{78} (2001),
  no.~1-2, 153--192.

\bibitem{Kilbas-Saigo-95}
A.A.~Kilbas and M.~Saigo, \emph{{On solution of integral equations of
  Abel-Volterra type}}, Differential Integral equations \textbf{8} (1995),
  no.~5, 993--1011.

\bibitem{Lamperti-72}
J.~Lamperti, \emph{semi-stable {M}arkov processes}, Z. Wahrsch. Verw. Geb.
  \textbf{22} (1972), 205--225.

\bibitem{Lebedev-72}
N.N.~Lebedev, \emph{Special functions and their applications}, Dover
  Publications, New York, 1972.

\bibitem{Levy-50}
P.~L{\'e}vy, \emph{Wiener's random functions, and other {L}aplacian random
  functions}, Proc. Sec. Berkeley Symp. Math. Statist. Probab., Vol. {II},
  Univ. California Press, Berkeley, 1950, pp.~171--186.

\bibitem{Maulik-Zwart-06}
K.~Maulik and B.~Zwart, \emph{Tail asymptotics for exponential functionals of
  {L}\'evy processes}, Stochastic Process. Appl. \textbf{116} (2006), 156--177.

\bibitem{Meyer-69}
P.A.~Meyer,
\emph{Processus \`a accroissements ind\'ependants et positifs}, S\'eminaire de probabilit\'es de Strasbourg \textbf{3} (1969),  175--189.


\bibitem{Mittag-03}
G.~Mittag-Leffler, \emph{Sur la nouvelle function ${E}_{\alpha}(x)$}, C. R.
  Math. Acad. Sci. Paris \textbf{137} (1903), 554--558.

\bibitem{Patie-06-poch}
P.~Patie, \emph{Law of the exponential functional of a family of one-sided {L}\'evy processes via self-similar continuous state branching processes with immigration and the {W}right hypergeometric functions}, Submitted (2007).

\bibitem{Pitman-Yor-80}
J.~Pitman and M.~Yor, \emph{{{B}essel processes and infinitely divisible
  laws}}, Stochastic Integrals (Proc. Sympos. Univ. Durham, Durham, 1980),
  Lecture Notes in Math. (D.~Williams, ed.), vol. 851, Springer, Berlin, 1981,
  pp.~285--370.

\bibitem{Revuz-Yor-99}
D.~Revuz and M.~Yor, \emph{Continuous martingales and brownian motion},
  $3^{rd}$ ed., vol. 293, Springer-Verlag, Berlin-Heidelberg, 1999.

\bibitem{Rivero-05}
V.~Rivero, \emph{Recurrent extensions of self-similar {M}arkov processes and
  {C}ram\'er's condition}, Bernoulli \textbf{11} (2005), no.~3, 471--509.

\bibitem{Sato-99}
K.~Sato, \emph{L\'evy processes and infinitely divisible distributions},
  Cambridge University Press, Cambridge, 1999.

\bibitem{Steutel-vanHarn-04}
F.W.~Steutel  and K.~van Harn,  \emph{Infinite divisibility of probability distributions on the real line},
  Monographs and Textbooks in Pure and Applied Mathematics, vol. 259, Marcel Dekker Inc., New York, 2004.


\bibitem{Vuolle-Apiala-94}
J.~Vuolle-Apiala, \emph{It\^o excursion theory for self-similar markov
  processes}, Ann. Probab., Vol. 22, No. 2 (1994), 546--565.

\bibitem{Wolfe-71}
S.J.~Wolfe, \emph{On the unimodality of {L} functions}, Ann. Math. Stat.
  \textbf{42} (1971), 912--918.

\bibitem{Wolfe-82}
S.J.~Wolfe, \emph{On a continuous analogue of the stochastic difference equation
  $\textrm{X}\sb{n}=\rho \textrm{X}\sb{n-1}+\textrm{B}\sb{n}$},
Stochastic Process. Appl., \textbf{12} (1982), no.~3 , 301--312.

\bibitem{Yor-80}
M.~Yor, \emph{Loi de l'indice du lacet brownien et distribution de
  {Hartman-Watson}}, Z. Wahrsch. Verw. Geb. \textbf{53} (1980), 71--95.

\bibitem{Yor-01}
\bysame, \emph{{Exponential functionals of Brownian motion and related
  processes}}, Springer Finance, Berlin, 2001.

\end{thebibliography}
\end{document}